\documentclass[reqno,12pt]{amsart}
\usepackage[latin1]{inputenc}
\usepackage{amsmath,amsthm,amssymb,amscd,pb-diagram}
\advance\evensidemargin -.8700in

\newcommand{\comment}[1]{}

\newcommand{\field}[1]{\mathbb{#1}}
\newcommand{\Q}{\field{Q}}

\newcommand{\F}{\field{F}}

\newtheorem{theorem}{Theorem}
\newtheorem{lemma}[theorem]{Lemma}
\newtheorem{prop}[theorem]{Proposition}

\theoremstyle{remark}

\theoremstyle{definition}

\newtheorem*{rem*}{Remark}
\newtheorem*{rems*}{Remarks}

\newtheoremstyle{citing}
  {3pt}
  {3pt}
  {\itshape}
  {}
  {\bfseries}
  {.}
  {.5em}
  {\thmnote{#3}}
\theoremstyle{citing}

\numberwithin{theorem}{section}
\numberwithin{equation}{section}

\begin{document}
\author{Daniel Rabayev}
\address{
Department of Mathematics\\
Technion --- Israel Institute of Technology\\
Haifa, 32000\\
Israel }
\email{danielr@tx.technion.ac.il}
\author{Jack Sonn}
\address{
Department of Mathematics\\
Technion --- Israel Institute of Technology\\
Haifa, 32000\\
Israel }
\email{sonn@tx.technion.ac.il}

\title[2-coverable symmetric and alternating groups]{On Galois realizations of the 2-coverable symmetric and alternating groups}
\date{\today}

\begin{abstract}
Let $f(x)$ be a monic polynomial in $\mathbb{Z}[x]$ with no
rational roots but with roots in $\mathbb{Q}_p$ for all $p$, or
equivalently, with roots mod $n$ for all $n$. It is known that
$f(x)$ cannot be irreducible but can be a product of two or more
irreducible polynomials, and that if $f(x)$ is a product of $m>1$
irreducible polynomials, then its Galois group must be ``$m$-coverable", i.e. a union of
conjugates of $m$ proper subgroups, whose total intersection is trivial.  We are thus led to a variant of the inverse Galois problem: given an $m$-coverable finite group $G$, find a Galois realization of $G$ over the rationals $\mathbb{Q}$ by a polynomial $f(x) \in \mathbb{Z}[x]$ which is a product of $m$ nonlinear irreducible factors (in $\mathbb{Q}[x]$) such that $f(x)$ has a root in $\mathbb{Q}_p$ for all $p$.
The minimal value $m=2$ is
 of special interest.  It is known that the symmetric group $S_n$ is $2$-coverable if and only if $3\leq n\leq6$, and the alternating group $A_n$ is $2$-coverable if and only if $4\leq n\leq8$.
 In this paper we solve the above variant of the inverse Galois problem for the $2$-coverable symmetric and alternating groups, and exhibit an explicit polynomial for each group, with the help of the software packages MAGMA, PARI and GAP.  

\end{abstract}
\maketitle
\section{Introduction}
Let $f(x)$ be a monic polynomial in $\mathbb{Z}[x]$ with no
rational roots but with roots in $\mathbb{Q}_p$ for all $p$, or
equivalently, with roots mod $n$ for all $n$. It is known that
$f(x)$ cannot be irreducible but can be a product of two or more
irreducible polynomials, and that if $f(x)$ is a product of $m>1$
irreducible polynomials, then its Galois group must be $m$-{\it coverable}, i.e. a union of
conjugates of $m$ proper subgroups, whose total intersection is trivial.  We are thus led to a variant of the inverse Galois problem: given an $m$-coverable finite group $G$, find a Galois realization of $G$ over the rationals $\mathbb{Q}$ by a polynomial $f(x) \in \mathbb{Z}[x]$ which is a product of $m$ nonlinear irreducible factors (in $\mathbb{Q}[x]$) such that $f(x)$ has a root in $\mathbb{Q}_p$ for all $p$.
\smallskip

 The existence of such a polynomial can be established by means of the following proposition.

\begin{prop}\label{prop:criterion} \cite{sonn} Let $K/\mathbb{Q}$ be a finite Galois extension
with Galois group $G$. The following are equivalent:

$ (1)$ $K$ is a splitting field of a product $f(x)=g_{1}(x)\cdot\cdot\cdot g_{m}(x)$
of $m$ irreducible polynomials of degree greater than one in $\mathbb{Q}[x]$
and $f(x)$ has a root in $\mathbb{Q}_{p}$ for all primes $p$.

$ (2)$ $G$ is the union of conjugates of $m$ proper
subgroups $A_{1},...,A_{m}$, the intersection of all these conjugates
is trivial, and for all primes $\mathfrak{p}$ of $K$, the decomposition
group $G(\mathfrak{p})$ is contained in a conjugate of some $A_{i}$.
\end{prop}

\begin{rem*} Condition (2) is evidently satisfied if $G$ is $m$-coverable and all its decomposition groups are cyclic.  This last condition holds automatically at all unramified primes. \end{rem*}

This remark is used in \cite{sonn} to prove existence for any $m$-coverable finite solvable group.

The minimal value $m=2$ is of natural interest.  All Frobenius groups are $2$-coverable.  In \cite{sonn} it is proved that Galois realizations satisfying condition (2) of Prop. 1.1 exist for all nonsolvable Frobenius groups with $m=2$.

In this paper we focus on $2$-coverable  symmetric and alternating groups, which have been determined in \cite{BBH} and \cite{Bubboloni}:  $S_n$ is $2$-coverable if and only if $3\leq  n \leq 6$, and $A_n$ is $2$-coverable if and only if $4\leq  n \leq 8$.

Polynomials for $S_3$ are exhibited in \cite{berend-bilu}.  For $S_4$ and $A_4$ existence follows from the general theorem for solvable groups in \cite{sonn}.  For $S_5$, existence follows from a known quintic (see e.g. \cite{sonnfrob}).  Until now, not even existence was known for $S_6$.  In the present paper, with the help of the software packages MAGMA, GAP and PARI, we both prove existence and produce an explicit polynomial for each $S_n$ ($3\leq  n \leq 6$) and each $A_n$ ($4\leq  n \leq 8$):

\begin{theorem}\label{main} Let $G$ be a $2$-coverable group which is either a symmetric group $S_n$ (i.e. $3\leq  n \leq 6$), or an alternating group $A_n$  (i.e. $4\leq  n \leq 8$).  Then $G$ is realizable over $\Q$ as the Galois group of a polynomial $f(x)$ which is the product of two (nonlinear) irreducible polynomials in $\Q[x]$, such that $f(x)$ has a root in $\Q_p$ for all primes $p$.
\end{theorem}

In each of the following sections, we take a group  $G \in \{S_4,S_5,S_6, A_4,...,A_8 \}$, and a $2$-covering of $G$ from \cite{BBH} or \cite{Bubboloni}, given by two subgroups $U_1,U_2$ of $G$.  We then produce a polynomial $f(x)\in \Q[x]$ with Galois group $G$ over $\Q$ which satisfies the condition that for
every prime $\mathfrak{p}$ of the splitting field $K$ of $f(x)$ over $\Q$, the decomposition
group $G(\mathfrak{p})$ is contained in a conjugate of $U_1$ or $U_2$.  We then find two irreducible polynomials $g_1(x),g_2(x)\in \Q[x]$ whose splitting fields are contained in $K$, such that for each $i=1,2$, if $\eta_i\in K$ is a root of $g_i(x)$,   then the Galois group of $K$ over $\Q(\eta_i)$ is $U_i$.  It then follows that $h(x)=g_1(x)g_2(x)$ has the desired property.  The software packages MAGMA, GAP and PARI are used in the searches for $f(x)$ and in the computation of $g_1(x),g_2(x)$.

\section{The symmetric group $S_{4}$}

The 2-covering of $S_{4}$ in \cite{BBH} is
given by $U_{1}$:= the stabilizer in $S_{4}$ of one
point and $U_{2}$:= the 2-Sylow subgroup of $S_{4}$.  $U_{1}$ contains a permutations of type%
\footnote{If $a_{1},...,a_{k}$ are natural numbers, and $\sigma$ is a permutation, then $\sigma$ is of type $[a_{1},...,a_{k}]$ if and only if $\sigma$ is the product of $k$ disjoint cycles
of length $a_{1},...,a_{k}$ respectively. In particular, a permutation of type $[n]$ is an $n$-cycle.%
} $[2]$ and  of type $[3]$, and $U_{2}$ contains  permutations of
types $[4]$ and $[2,2]$. Since two elements of $S_{n}$
are conjugate if and only if they are of the same type, $U_1,U_2$ is a 2-covering of $S_4$.

Now that we have the covering we can introduce the polynomial
\begin{eqnarray*}
f(x) & = & x^{4}-5x^{2}+x+4\\
\end{eqnarray*}
 with prime discriminant  $p=2777$.  $f(x)$ is irreducible mod
$3$, factors into the product $(x+4)(x^{3}+7x^{2}+1)$ of irreducible
factors mod $11$ and factors into the product $(x+8)(x+19)(x^{2}+19x+20)$
mod $23$. These factorizations imply that the Galois group of $f(x)$
contains a 4-cycle, 3-cycle and a transposition. This implies that
the Galois group of $f(x)$ (over $\Q$), which is contained in $S_{4}$, is indeed
$S_{4}$.

Let $K$ be the splitting field of $f(x)$ over $\Q$. We will show that for all primes $\mathfrak{p}$ of $K$, the
decomposition group $G(\mathfrak{p})$
is cyclic, hence contained in a conjugate of $U_{1}$ or $U_{2}$.
The decomposition group of an unramified prime is cyclic, so it
is enough to check the ramified primes. If the
discriminant of the polynomial is squarefree then the ramified primes
are exactly those primes dividing the discriminant of the polynomial, which in our case is a prime $p=2777$.

We shall show that $G(\mathfrak{p})$ is cyclic  by showing that the Galois group of $f(x)$ over $\mathbb{Q}_{p}$, which is isomorphic to $G(\mathfrak{p})$,
is of order 2. We have the following factorization of $f(x)$ modulo  $p=2777$: \begin{eqnarray*}
 & f(x)\equiv(x+787)(x+929)(x+1919)^{2}\ (mod\ 2777)\\
\end{eqnarray*}
Applying Hensel's Lemma we conclude that $f(x)$ factors into two
relatively prime linear factors and one irreducible factor of degree
2 over $\mathbb{Q}_{p}$. Indeed, Hensel's Lemma tells us that $f(x)$
will factor into two relatively prime linear factors and one factor
of degree 2, which must be irreducible since
$f(x)$ is separable over $\mathbb{Q}_p$ and $p$ is ramified. This implies that the Galois
group of $f(x)$ over $\mathbb{Q}_{p}$ is of order 2.

We now find two irreducible polynomials $g_{1}(x),g_{2}(x)$
of degree greater than 1, such that
$K$ is the splitting field of their product and the groups $U_{1}$
and $U_{2}$ are the Galois groups of $K$ over the root fields of $g_{1}(x)$
and $g_{2}(x)$ respectively. Moreover, the polynomial $g_{1}(x)g_{2}(x)$
has a root in $\mathbb{Q}_{p}$ for all prime $p$.

Since $U_{1}$ is the stabilizer of one point, we may take $g_1(x):=f(x)$.

Since $U_{2}$ is a 2-sylow subgroup of $S_4$, we may take $g_{2}(x)$ to be
the cubic resolvent \begin{eqnarray*}
g_{2}(x) & = & x^{3}+10x^{2}+9x+1\\
\end{eqnarray*} of $f(x)$.

In conclusion, the product:\begin{eqnarray*}
g_{1}(x)\cdot g_{2}(x) & = & (x^{4}-5x^{2}+x+4)(x^{3}+10x^{2}+9x+1)\\
\end{eqnarray*}
has Galois group $S_{4}$ and has a root in $\mathbb{Q}_{p}$ for
all prime $p$.

\section{The symmetric group $S_{5}$}


The 2-covering of $S_{5}$ in \cite{BBH} is given
by $U_{1}$ := the stabilizer in $S_{5}$ of the set \{1,2\},
and $U_{2}:=$ the normalizer $N_{S_{5}}\{\langle(1,2,3,4,5)\rangle\}$ of $\langle(1,2,3,4,5)\rangle$ in $S_5$.  $U_{1}$
contains permutations of type $[2]$, $[3]$ and $[2;3]$.
$U_{2}$ contains  permutations of type $[4]$, $[5]$ and $[2;2]$.
Therefore $U_1,U_2$ is a 2-covering of $S_5$.

Let \begin{eqnarray*}
f(x) & = & x^{5}-9x^{3}-9x^{2}+3x+1,\\
\end{eqnarray*}
a polynomial with prime discriminant $p=36497$.
 $f(x)$ is irreducible mod
$17$, factors into the product $(x+1)(x^{4}+2x^{3}+x^{2}+2x+1)$
of irreducible factors mod $3$ and factors into the product $(x+37)(x+56)(x+76)(x^{2}+37x+48)$
mod $103$. These factorizations imply that the Galois group of $f(x)$
contains a 5-cycle, 4-cycle and a transposition, and is therefore
$S_{5}$.

As in the previous section we will show
that for all primes $\mathfrak{p}$ of the splitting field $K$ of $f(x)$, the decomposition group
$G(\mathfrak{p})$ is cyclic and thus contained in a conjugate of
$U_{1}$ or $U_{2}$. Again it is enough to check the ramified primes.   The discriminant of $f(x)$ is a prime $p=36497$, so there is only one ramified prime.

 The factorization of $f(x)$ modulo $p$ is  \begin{eqnarray*}
 & f(x)\equiv(x+8522)(x+18488)(x+19525)(x+31478)^{2}\ \ (mod\ 36497).\\
\end{eqnarray*}
Applying Hensel's Lemma as in the preceding section, the decomposition group at $p$ is of order 2, hence cyclic.

We now find the desired polynomials $g_1(x),g_2(x)$. We  start
with the polynomial $g_1(x)$ corresponding to the group $U_{1}$, which is the stabilizer in $S_{5}$ of the set \{1,2\}.
It is easy to see that $U_{1}=C_{2}\times S_{3}$, and furthermore, all
of the subgroups
 of $S_{5}$ which are isomorphic to $C_{2}\times S_{3}$ are conjugates.
Hence we need only  find a polynomial whose roots lie in $K$ and with Galois group $U_{1}$
over its root field. That polynomial is the polynomial whose roots
are the sums of pairs of roots of $f(x)$, namely (using GAP), \begin{eqnarray*}
g_{1}(x) & = & x^{10}-27x^{8}-9x^{7}+234x^{6}+151x^{5}-756x^{4}-585x^{3}+873x^{2}+660x-163\\
\end{eqnarray*}
Indeed, it is an irreducible polynomial of degree 10 ($=[S_5:U_1]$) and according
to the previous discussion, it is clear that it has Galois group $C_{2}\times S_{3}$
over its root field.\\
\\
We turn next to the polynomial $g_2(x)$ corresponding to the group $U_{2}=N_{S_{5}}\langle(1,2,3,4,5)\rangle$.
The order of $U_{2}$ is 20, so $g_2(x)$ will be of degree 6.  First, let us notice that up
to conjugation, there is only one subgroup of order 20 in $S_{5}$,
since every subgroup of order 20
is the normalizer of a 5-cycle.   Thus we need only find an irreducible polynomial
of degree 6 whose splitting field is contained in $K$.  Our method is ad hoc.
Consider the action of $S_5$ on the set $A$ of the ten 2-subsets of $\{1,2,3,4,5\}$, which is equivalent to its
action on the roots of $g_1(x)$.  Then consider the action of $S_5$ on the set $B$ of 252 5-subsets of $A$.  This latter action
has an orbit of length 12 (MAGMA).  
We wish to produce a separable polynomial of degree 252 with splitting field $K$ such that the action of $G(K/\Q)\cong S_5$ on its roots is equivalent to its action on $B$.
We start with
 \begin{eqnarray*}
g_1(x) & = & x^{10}-27x^{8}-9x^{7}+234x^{6}+151x^{5}-756x^{4}-585x^{3}+873x^{2}+660x-163\\
\end{eqnarray*}
as above which is an irreducible polynomial with Galois group $S_{5}$ over $\Q$.
Using a Tschirnausen algorithm%
\footnote{We are grateful to Alexander Hulpke for introducing us to this algorithm in GAP.%
} in GAP, we obtain a degree 10 polynomial $g(x)$ the sum of whose roots is zero and with the same root field (up to conjugacy) as $g_1(x)$ and such that
the polynomial $h(x)$ of degree 252 whose roots are the sums of quintuples of distinct roots of $g(x)$, is separable.  We now observe that since the sum of the roots of $g(x)$ is zero, the negative of any root of $h(x)$ is also a root of $h(x)$, so $h(x)=r(x^2)$ for some $r(t)\in\Q[x]$ of degree 126.  $r(t)$ has an irreducible factor of degree
6 with splitting field contained in $K$, which is
$$g_{2}(x) =  x^{6}-\frac{1389834}{26569}x^{5}+\frac{804440497905}{705911761}x^{4}
  -\frac{248200341684305820}{18755369578009}x^{3} $$

   $$ +\frac{43053809146811239035780}{498311414318121121}x^{2}
  -\frac{3981047272192561628814509657}{13239635967018160063849}x$$

$$  +\frac{153302536976991922308475759105476}{351763888007705494736404081}$$

and is therefore a polynomial of the desired type.

\section{The symmetric group $S_{6}$}

The 2-covering of $S_{6}$ in \cite{BBH} is given
by $U_{1}$ := the stabilizer in $S_{6}$ of the
partition \{1,2,3\},\{4,5,6\} and $U_{2}=S_{5}$.   $U_{1}$ contains
a permutation of type $[6]$, hence also a permutation of type $[2,2,2]$
and $[3,3]$. Moreover, $U_{1}$ also contains a permutation of type
$[2,4]$, thus for each type of fixed-point-free permutation, $U_{1}$
contains at least one representative. On the other hand, $U_{2}$
contains representatives for each type of permutation with at least
one fixed point. Therefore $U_1,U_2$ is a 2-covering of $S_6$.  Note that $U_1$ is of order 72 and is also the normalizer of a 3-Sylow
subgroup of $S_6$.

 The polynomial
\begin{eqnarray*}
f(x) & = & x^{6}-10x^{4}-9x^{3}+5x^{2}+3x-1\\
\end{eqnarray*}
has prime discriminant $p=33994921$,
 is irreducible mod
$37$, factors into the product $(x+11)(x^{5}+2x^{4}+7x^{3}+5x^{2}+2x+7)$
of irreducible factors mod $13$ and factors into the product $(x+11)(x+70)(x+199)(x+232)(x^{2}+14x+130)$
mod $263$. It follows that the Galois group of $f(x)$
contains a 6-cycle, 5-cycle and a transposition, so its
Galois group   is
$S_{6}$.

Next we will show that for all primes $\mathfrak{p}$ of the splitting field $K$ of $f(x)$, the
decomposition group $G(\mathfrak{p})$ is cyclic and thus contained
in a conjugate of $U_{1}$ or $U_{2}$. Again it is enough to check the ramified
primes. Again the discriminant is a prime $p=33994921$,  so there is only one ramified prime.

We shall now show that the Galois group of $f(x)$ over $\mathbb{Q}_{p}$
is cyclic of order 2, hence so is $G(\mathfrak{p})$. $f(x)$ factors into a product of four distinct linear
factors and the square of a linear factor  modulo $p$; explicitly,
 \begin{eqnarray*}
f(x)\equiv & (x+665896)(x+3641312)(x+15713959)(x+25142575)(x+11413050)^{2}\ \ (mod\ p)\end{eqnarray*}
Applying Hensel's Lemma as in the preceding sections, the Galois group of $f(x)$ over $\mathbb{Q}_{p}$
is of order 2 and thus contained in a conjugate of $U_{1}$ or $U_{2}$.

We will now find the polynomials $g_{1}(x)$ and $g_{2}(x)$. Let
us notice that we already have a polynomial for one of the conjugates
of the group $U_{2}$, namely $g_{2}(x):=f(x)=x^{6}-10x^{4}-9x^{3}+5x^{2}+3x-1$,
 since the Galois group of $f(x)$ over $\mathbb{Q}(\theta)$,
where $\theta$ is a root of $f(x)$, is $S_{5}$.

We now turn to the polynomial $g_{1}(x)$ which corresponds to the group $U_{1}$,
 which
is the stabilizer of the partition $\{1,2,3\},\{4,5,6\}$.
Applying GAP we
construct a polynomial $h(x)$ whose roots are the sums of triples of distinct roots of $f(x)$:
\begin{eqnarray*}
h(x) & = & x^{20}-60x^{18}+1470x^{16}-19209x^{14}+146153x^{12}-662097x^{10}\\
 &  & +1751524x^{8}-2532942x^{6}+1684385x^{4}-278413x^{2}+4096.\end{eqnarray*}
Observe that $h(x)=g_1(x^2)$, where
  \begin{eqnarray*}
g_{1}(x) & = & x^{10}-60x^{9}+1470x^{8}-19209x^{7}+146153x^{6}-662097x^{5}\\
 &  & +1751524x^{4}-2532942x^{3}+1684385x^{2}-278413x+4096.\end{eqnarray*}
 $g_{1}(x)$ is irreducible over $\mathbb{Q}$.
It is evident that the splitting field of $h(x)$ is contained in
$K$, and since $g_{1}(x)$ is irreducible over $\Q$, the splitting field of $h(x)$ equals $K$.
It is evident that $h(x)$
 has Galois group $S_{3}\times S_{3}$
over $\mathbb{Q}(\eta)$, where $\eta$ is a root of $h(x)$.
The root field $\mathbb{Q}(\eta^2)$
of $g_{1}(x)$ is a subfield of $\mathbb{Q}(\eta)$ of degree $10$ over $\mathbb{Q}$
, so $\Gamma:=G(K/\mathbb{Q}(\eta^2))$ is a subgroup of order $72$ of $S_{6}$, containing $S_{3}\times S_{3}$.
$\Gamma$ contains a 3-Sylow subgroup of $S_6$ whose normalizer is of index 1 or 2 in $\Gamma$.  It now follows from the Sylow theorems the index must be 1, so $\Gamma$ must be conjugate to $U_1$ in $S_6$.  (In fact, every subgroup of $S_6$ of order 72 is conjugate to $U_1$.)

\section{The alternating group $A_{4}$}

A 2-covering of $A_{4}$ is given by $U_{1}$:=
the stabilizer in $A_{4}$ of one point, namely, $U_{1}=C_{3}$, and
the group $U_{2}$:=the 2-sylow subgroup of $A_{4}$ which
is $C_{2}\times C_{2}$. It is evident that this is indeed a 2-covering
of $A_{4}$.

Let \begin{eqnarray*}
f(x) & = & x^{4}-10x^{3}-7x^{2}+3x+2.\\
\end{eqnarray*}
The discriminant
of $f(x)$ is $163^{2}$, a square (of a prime), hence the Galois
group is contained in $A_{4}$.
 $f(x)$ factors into the product
$(x+1)(x^{3}+x^{2}+x+2)$ of irreducible factors mod $3$ and factors
into the product $(x^{2}+2x+4)(x^{2}+3x+3)$ mod $5$. These factorizations
imply that the Galois group of $f(x)$ contains a 3-cycle and a product
of two transpositions, hence the Galois group is  $A_{4}$.

Now we  check  that for all primes
$\mathfrak{p}$ of the splitting field $K$, the decomposition group $G(\mathfrak{p})$,
is cyclic and thus contained in a conjugate of $U_{1}$ or $U_{2}$.
As noted before,  it is enough to check the ramified primes. Since the discriminant is $p^{2}=163^{2}$, which is the square of a prime, the only prime to check is $p=163$.

To this end we shall show that the
Galois group of $f(x)$ over $\mathbb{Q}_{p}$ is of order 2.   $f(x)$ factors mod $p$ as: \begin{eqnarray*}
 & f(x)\equiv(x+50)(x+143)^{3}\ \ (mod\ 163)\\
\end{eqnarray*}
Applying Hensel's Lemma, $f(x)$ factors into a linear
factor and a factor of degree 3 over $\mathbb{Q}_{p}$, which we claim must be irreducible over $\mathbb{Q}_{p}$.  Indeed,  $f(x)$ is
separable over $\mathbb{Q}_{p}$  so this cubic factor is either irreducible
as claimed or a product of a linear factor and an irreducible factor
of degree 2. In the latter case,  $f(x)$ has two relatively prime linear factors which implies that the local Galois group fixes exactly two roots, thus contains a transposition which contradicts
the fact that this Galois group is contained in $A_{4}$. This proves the claim and implies
that the Galois group of $f(x)$ over $\mathbb{Q}_{p}$ is either
$S_{3}$ or $C_{3}$. The decomposition group cannot be $S_{3}$ since
it is a subgroup of $A_{4}$ and $A_{4}$ does not contain a copy
of $S_{3}$, hence we deduce that it is (the cyclic group) $C_{3}$.

We now proceed to the desired polynomials $g_1(x),g_2(x)$.  This is similar to the case of $S_4$.
Since $U_{1}$ is the stabilizer of one point, we may take $g_1(x):=f(x)$.
Since $U_{2}$ is a 2-sylow subgroup of $A_4$, we may take $g_{2}(x)$ to be
the cubic resolvent \begin{eqnarray*}
g_{2}(x) & = & x^{3}+89x^{2}+2586x+24649\\
\end{eqnarray*} of $f(x)$.

\section{The alternating group $A_{5}$}

A 2-covering of $A_{5}$ can be obtained from a 2-covering
of $S_{5}$ and the following easy lemma.

\begin{lemma} \cite{Bubboloni}. Let $G$ be a $2$-coverable group which is covered by
the conjugates of the subgroups $H$ and $K$. If $N\trianglelefteq G$ and $G=NH=NK$, then
$N$ admits the covering $H\cap N,\ K\cap N$.
\end{lemma}

Applying the lemma to the 2-covering of $S_{5}$ which was given in
a previous section, we take $U_{1}$ to be the intersection of
$A_{5}$ with the stabilizer $Stab_{S_{5}}\{1,2\}$  of the set \{1,2\} in $S_{5}$,
and $U_{2}=N_{S_{5}}\{\langle(1,2,3,4,5)\rangle\}\cap A_{5}$.
As $S_{5}=N_{S_{5}}\{\langle(1,2,3,4,5)\rangle\}\cdot A_{5}=Stab_{S_{5}}\{1,2\}\cdot A_{5}$,
$\{U_1,U_2\}$ is a 2-covering of $A_{6}$.

 Let \begin{eqnarray*}
f(x) & = & x^{5}-5x^{4}-7x^{3}+18x^{2}-x-1.\\
\end{eqnarray*}
The discriminant of
$f(x)$ is $15733^{2}$,  a square (of a prime), so the Galois group is contained in $A_{5}$.
 $f(x)$ is irreducible mod $3$, factors into the product
$(x+3)(x+4)(x^{3}+3x^{2}+2)$ of irreducible factors mod $5$ and
factors into the product $(x+5)(x^{2}+x+4)(x^{2}+3x+1)$ mod $7$.
These factorizations imply that the Galois group of $f(x)$ contains
a 5-cycle, a 3-cycle and a product of two transpositions, thus $2\cdot3\cdot5=30$
divides the order of the group. Since $A_{5}$ has no subgroup
of order 30, the Galois group of $f(x)$ is $A_{5}$.

Next we will show that for all primes $\mathfrak{p}$ of the splitting field $K$ of $f(x)$, the
decomposition group $G(\mathfrak{p})$ is cyclic and thus contained
in a conjugate of $U_{1}$ or $U_{2}$.  As usual it is enough to check the ramified
primes, and as the discriminant $15733^{2}=p^{2}$
 is the square of a prime, it suffices to check that $G(\mathfrak{p})$ is cyclic at this single $p$.

 Now \begin{eqnarray*}
 & f(x)\equiv_{}(x+13759)(x+4886)(x+4272)^{3}\ \ (mod\ 15733).\\
\end{eqnarray*}
The proof proceeds as in the preceding case.  By Hensel's Lemma, $f(x)$ factors into two
relatively prime linear factors and a factor of degree 3 over $\mathbb{Q}_{p}$, which we claim must be irreducible over $\mathbb{Q}_{p}$. Indeed,  $f(x)$ is
separable over $\mathbb{Q}_{p}$  so this cubic factor is either irreducible
as claimed or a product of a linear factor and an irreducible factor
of degree 2.  In the latter case,  $f(x)$ has 3 relatively prime linear factors, which implies that the local Galois group fixes exactly 3 roots, thus contains a transposition, which contradicts
the fact that this Galois group is contained in $A_{5}$. This proves the claim and implies
that the Galois group of $f(x)$ over $\mathbb{Q}_{p}$ is either
$S_{3}$ or $C_{3}$. Now the fact that there
are two linear factors means that the Galois group acts trivially
on two roots, thus the decomposition group is a subgroup of $A_{3}$.
But this implies that the Galois group of $f(x)$ over $\mathbb{Q}_{p}$
is $C_{3}$.

We turn now to the polynomials $g_1(x),g_2(x)$. From here the proof is similar to the case of $S_5$.
The indices of $U_{1}$ and $U_{2}$ in $A_{5}$ are 10 and 6 respectively.  Furthermore, each $U_i$ is the unique subgroup of order $|U_i|$ in $A_5$ up to conjugacy, $i=1,2$.  By the same method as in the case of $S_5$, starting with $f(x)$, we compute
\begin{eqnarray*}
g_1(x) & = & x^{10}-20x^{9}+129x^{8}-167x^{7}-1160x^{6}+3139x^{5}\\
 &  & +2214x^{4}-9559x^{3}+1089x^{2}+5486x-461\end{eqnarray*}
and
\begin{eqnarray*}
 & g_{2}(x)= & x^{6}-\frac{42323896}{212521}x^{5}+\frac{733217461235730}{45165175441}x^{4}-\frac{6645087421265775407450}{9598548249896761}x^{3}\\
 &  & +\frac{33179389816151522298285969465}{2039893072616309544481}x^{2}-\frac{86429398882716670964952901401177395}{433520115685490720702646601}x\\
 &  & +\frac{91715516822801665361140634800285822768801}{92132128505596173454447158291121}.\end{eqnarray*}

\section{The alternating group $A_{6}$}

In contrast to the previous examples, in this example not all the decomposition
groups turn out to be cyclic.
The 2-covering of $A_{6}$ is given by applying the lemma
from the previous section to the 2-covering of $S_{6}$ which was
given  earlier.  Namely, $U_{1}=H_{1}\cap A_{6}$, where
$H_{1}$ denotes the stabilizer in $S_{6}$ of the partition $\{1,2,3\}$,
$\{4,5,6\}$ and $U_{2}=H_{2}\cap A_{6}=A_{5}$ where $H_{2}=S_{5}$.
($S_{6}=A_{6}H_{1}=A_{6}H_{2}$ because $U_{1}$
and $U_{2}$ each contains an odd permutation.)

 Let \begin{eqnarray*}
f(x) & = & x^{6}-3x^{5}-4x^{4}+5x^{3}+3x^{2}-5x-1.\\
\end{eqnarray*}
The discriminant
of $f(x)$ is $14341^{2}$  is a square (of a prime), so the Galois group is contained in $A_{6}$.
 $f(x)$ factors as $(x^{2}+x+2)(x^{4}+2x^{3}+x^{2}+1)$
 mod $3$, as $(x+3)(x^{5}+4x^{4}+4x^{3}+3x^{2}+4x+3)$
mod $5$,  as $(x+16)(x+53)(x+61)(x^{3}+69x^{2}+8x+82)$
 mod $101$ and as $(x^{3}+13x^{2}+27x+47)(x^{3}+93x^{2}+68x+51)$
mod 109, hence the Galois group has order divisible by $9\cdot4\cdot5=180$.
Since $A_{6}$ has no subgroup of order 180, the Galois group is $A_{6}$.

In order to verify that the requirements hold, we will show that for
all primes $\mathfrak{p}$ of the splitting field $K$ of $f(x)$, the decomposition group $G(\mathfrak{p})$,
is contained in a conjugate of $U_{1}$ or $U_{2}$. Again it is enough to check this at
the (single) ramified prime $p=14341$.

We claim that the Galois group of $f(x)$ over
$\mathbb{Q}_{p}$ is contained in a conjugate of $U_{2}$.  $f(x)$ factors mod $p$ as a product of two distinct linear
factors and the square of an irreducible quadratic factor:  \begin{eqnarray*}
 & f(x)\equiv(x+5464)(x+5605)(x^{2}+8805x+9098)^{2}\ \ (mod\ 14341)\\
\end{eqnarray*}
By Hensel's Lemma, $f(x)$ factors into two
relatively prime linear factors and one irreducible factor of degree
4 over $\mathbb{Q}_{p}$. This means that the decomposition group,
which is contained in $A_{6}$, fixes two roots of $f(x)$, and is therefore contained in a conjugate of $A_5=U_2$.  Note that the decomposition group is noncyclic, since its order is divisible by $4$ and is embedded in $A_4$.

We turn to the polynomials $g_{1}(x)$ and $g_{2}(x)$.  First, as $U_2=A_5$, we may take $g_{2}(x)=f(x)=x^{6}-3x^{5}-4x^{4}+5x^{3}+3x^{2}-5x-1$.  By  the same method used for $S_6$, we arrive at
the polynomial
\begin{eqnarray*}
g_{1}(x) & = & x^{10}-\frac{93}{2}x^{9}+\frac{14253}{16}x^{8}-\frac{75007}{8}x^{7}+\frac{7698057}{128}x^{6}-\frac{62688663}{256}x^{5}+\frac{1276887561}{2048}x^{4}\\
 &  & -\frac{1886832143}{2048}x^{3}+\frac{45377978093}{65536}x^{2}-\frac{33044483389}{131072}x+\frac{1696039489}{1048576}\end{eqnarray*}
\section{The alternating group $A_{7}$}

The 2-covering of $A_7$ given in \cite{Bubboloni} is $U_1=N_{A_7}\langle(1234567)\rangle$, of order 21, and $U_2=(Sym\{1,2\}\times Sym\{3,4,5,6,7\})\cap A_7$, of order 120, where $Sym\{1,2\}$ is the symmetric group on $\{1,2\}$ and $Sym\{3,4,5,6,7\}$ is the symmetric group on $\{3,4,5,6,7\}$. Let \begin{eqnarray*}
f(x) & = & x^{7}-2x^{6}-5x^{5}-x^{4}-3x^{3}-x^{2}-x-5.\end{eqnarray*}
Since the discriminant
of the polynomial is $554293^{2}$, the square of a prime, the Galois group of $f(x)$ over $\mathbb{Q}$ is contained
in $A_{7}$.  The factorizations of $f(x)$ modulo $3,5,7,23$ show that the Galois group contains permutations of the types [7],[5],[2,4],[3,3] respectively.  As $A_7$ has no proper subgroup of order divisible by $3\cdot 4 \cdot 5 \cdot 7$ (GAP), the Galois group is $A_7$.

We show next that the Galois group of $f(x)$ over $\mathbb{Q}_{p}$
is contained in a conjugate of $U_{2}$, where $p=554293$ is the unique ramified prime. (As before, this is the only local Galois group that needs to be checked.)  The factorization of $f(x)$ modulo $p$ is  \begin{eqnarray*}
 & f(x)\equiv(x+521396)(x+134869)^{2}(x+281696)^{2}(x^{2}+308251x+256808)\ \ (mod\ 554293).\end{eqnarray*}
Applying Hensel's Lemma we conclude that $f(x)$ factors into a (separable) product
of  a linear factor and  quadratic factors
over $\mathbb{Q}_{p}$. It follows that the decomposition group $G(\mathfrak{p})$
is a 2-group. But this means that the decomposition group is contained
in a conjugate of $U_{2}$ since $U_{2}$ contains a 2-sylow subgroup
of $A_{7}$.

For the explicit polynomials $g_1(x),g_2(x)$, we start with $g_1(x)$ and observe first that $U_1$ is not a maximal subgroup of $A_7$.  It is contained in a maximal subgroup of order 168, isomorphic to $GL_3(2)$.  Although the unique conjugacy class of  subgroups of order 168 in $S_7$ splits into two conjugacy classes in $A_7$, each contains the unique conjugacy class of subgroups of order 21 in $A_7$.  It follows that either of the two conjugacy classes of subgroups of order 168 in $A_7$ can replace the conjugacy class of $U_1$ in the given 2-covering of $A_7$.  It is advantageous to do this in order to seek a $g_1(x)$ of smaller degree (15) than that of the $g_1(x)$ (120) that would correspond to the original $U_1$.  The method for seeking the explicit $g_1(x)$ and $g_2(x)$ is similar to the preceding, except for a greater dependence on trial and error, and that in addition to the two algorithms in GAP used earlier, the Tschirnhausen algorithm and the algorithm that  produces a polynomial whose roots are sums of roots of a given polynomial, an algorithm that  produces a polynomial whose roots are products of roots of a given polynomial was pressed into service.

The polynomials corresponding to the groups $U_{1}$ and $U_{2}$
are:

\begin{eqnarray*}
 & g_{1}(x)=x^{15}+\frac{236441}{835}x^{14}-\frac{553005}{3796}x^{13}-\frac{142958}{24613}x^{12}+\frac{3151297}{45345}x^{11}-\frac{5686331}{48267}x^{10}\\
 & +\frac{30589}{7924}x^{9}+\frac{10004127}{37216}x^{8}-\frac{255049}{586}x^{7}-\frac{41630027}{36825}x^{6}+\frac{72184440}{8779}x^{5}\\
 & -\frac{45952163}{35584}x^{4}-\frac{1964679}{5893}x^{3}+\frac{14351259}{5423}x^{2}-\frac{3135356}{2567}x-\frac{791507}{457}\end{eqnarray*}

\begin{eqnarray*}
 & g_{2}(x)=x^{21}-30x^{20}+395x^{19}-2937x^{18}+12917x^{17}-29077x^{16}-6543x^{15}\\
 & +236441x^{14}-553005x^{13}-142958x^{12}+3151297x^{11}-5686331x^{10}\\
 & +305889x^{9}+10004127x^{8}-255049x^{7}-41630022x^{6}+72184440x^{5}\\
 & -45952163x^{4}-1964679x^{3}+14351259x^{2}-3135356x-791507\end{eqnarray*}
respectively.

\section{The alternating group $A_{8}$}

The 2-covering of $A_8$ in \cite{Bubboloni} is given by $U_1$:=the affine linear group on $\F_2^3$ acting as permutations on the 8 points of this space (shown to be embedded into $A_8$), isomorphic to $\F_2^3\rtimes GL_3(\F_2)$, of index 15, and $U_2:=[Sym\{1,2,3\}\times Sym\{4,5,6,7,8\}]\cap A_{8}$, of index 56.
Let \begin{eqnarray*}
f(x) & = & x^{8}-x^{7}-2x^{6}-x^{5}+3x^{4}+3x^{3}+2x^{2}+x+1.\end{eqnarray*}
Since the discriminant
of the polynomial is $11489^{2}$, the square of a prime,  the Galois group of $f(x)$ over $\mathbb{Q}$ is contained
in $A_{8}$.  The factorizations of $f(x)$ modulo $3,37,41$ show that the Galois group contains permutations of the types $[7],[4,4],[3,5]$ respectively.  As $A_8$ has no proper subgroup of order divisible by $3\cdot 4 \cdot 5 \cdot 7$ (GAP), the Galois group is $A_8$.

Now
$f(x)$ factors modulo $p$ as: \begin{eqnarray*}
 & f(x)\equiv(x+1440)^{3}(x^{2}+10435x+8884)(x^{3}+8222x^{2}+9218x+10584)\ \ (mod\ 11489)\end{eqnarray*}

  Hensel's Lemma yields a corresponding factorization  $f(x)=a(x)b(x)c(x)$ over $\mathbb{Q}_p$ with $a(x),c(x)$ of degree 3 and $b(x)$ of degree 2.  The splitting field $K_p$ of $f(x)$ over $\mathbb{Q}_p$ is the composite of the splitting fields $K_a,K_b,K_c$ of these three factors, the first of which is ramified and the others unramified.  $a(x)$ cannot be a product of linear factors since $p$ ramifies in the splitting field of $f(x)$.  Thus $a(x)$ is either an irreducible cubic or a product of a linear and a quadratic polynomial.  $b(x)$ and $c(x)$ are irreducible.  If $a(x)$ is
 a product of a linear and a quadratic polynomial, then $K_a/\mathbb{Q}_p$ is a ramified quadratic extension and $G(K_p/\mathbb{Q}_p)\cong C_2\times C_2 \times C_3$, and contains transpositions, contrary to the fact that it is contained in $A_8$.  Therefore $a(x)$ is  an irreducible cubic.  It follows that $G(K_a/\Q_p)$ is $C_3$ or $S_3$.  $C_3$ is not possible because $K_a/\Q_p$ would then be totally and tamely ramified, but $p \equiv 2$ (mod 3) hence $\Q_p$ does not contain the cube roots of unity.  Hence  $G(K_a/\Q_p)$ is  $S_3$.  The only way that $G(K_p/\mathbb{Q}_p)$ can be contained in $A_8$ is if $K_b$ is contained in $K_a$ as the unramified quadratic subfield, and the action of $G(K_a/\Q_p)$ on the roots of $a(x)b(x)$ is the representation of $S_3$ as $\langle (123),(12)(45) \rangle$.  Furthermore, since $K_c/\Q_p$ is unramified of degree $3$, $K_a\cap K_c=\Q_p$, so $G(K_p/\mathbb{Q}_p)$ acts on the roots of $f(x)$ as $\langle(123),(12)(45),(678)\rangle$, which is contained in $U_2$.

 A trial and error search similar to that in the preceding section yields the following polynomials corresponding to the groups $U_{1}$ and $U_{2}$:

\begin{eqnarray*}
 & g_{1}(x)=x^{15}-\frac{95}{2}x^{14}+\frac{3197}{4}x^{13}-\frac{210405}{32}x^{12}-\frac{7433761}{256}x^{11}+\frac{22263473}{256}x^{10}-\frac{104363621}{512}x^{9}\\
 & -\frac{11710437647}{4096}x^{8}+\frac{918294939277}{65536}x^{7}+\frac{5512444210853}{131072}x^{6}-\frac{21829527485865}{262144}x^{5}-\frac{410948046145887}{2097152}x^{4}\\
 & -\frac{14737903234907309}{16777216}x^{3}-\frac{27794916917226379}{8388608}x^{2}-\frac{31865126119189869}{4194304}x-\frac{1256672780150364253}{134217728}\end{eqnarray*}

and

\begin{eqnarray*}
 & g_{2}(x)=x^{56}-21x^{55}+180x^{54}-735x^{53}+783x^{52}+5407x^{51}-21690x^{50}+6260x^{49}+134354x^{48}\\
 & -256270x^{47}-304036x^{46}+1275397x^{45}+942277x^{44}-6877457x^{43}-2069019x^{42}+42340751x^{41}\\
 & -37776712x^{40}-158592858x^{39}+379188834x^{38}+107601914x^{37}-1455434817x^{36}+1392744196x^{35}\\
 & +2854906257x^{34}-7385711268x^{33}+817821146x^{32}+17194579246x^{31}-19188201050x^{30}\\
 & -18578454420x^{29}+56645526134x^{28}-18919098068x^{27}-77381156815x^{26}+94740276783x^{25}\\
 & +30587068228x^{24}-150642640537x^{23}+94400959739x^{22}+77653585985x^{21}-155812844138x^{20}\\
 & +65943540240x^{19}+68936866008x^{18}-121634333784x^{17}+85946585062x^{16}-14385996295x^{15}\\
 & -52121609885x^{14}+76576460155x^{13}-37676127233x^{12}-34141244957x^{11}+82692304838x^{10}\\
 & -86059996214x^{9}+64840248670x^{8}-43929749632x^{7}+29390580877x^{6}-17890022300x^{5}\\
 & +8764351999x^{4}-3189669694x^{3}+811672328x^{2}-129075408x+9546832\end{eqnarray*}
respectively.

\bibliographystyle{plain}

\def\cprime{$'$}

\end{document}